\documentclass[10pt,a4paper,twoside]{article}

\title{\sc A Canonical Extension of Korn's First Inequality to {\sf H}{\rm(Curl)}
motivated by Gradient Plasticity\\ with Plastic Spin}
\def\shorttitle{Korn's inequality in {\sf H}{\rm(Curl)}}
\def\pauthor{Patrizio Neff, Dirk Pauly, Karl-Josef Witsch}

\def\mylabelonoff{off}
\def\allowdisbrk{no}

\usepackage{a4,exscale,ifthen,amsfonts,amssymb,amsmath,amscd,graphicx}
\usepackage[english]{babel}
\usepackage[mathscr]{eucal}

\author{{\sf\pauthor}}
\numberwithin{equation}{section}

\newcommand{\bewboxw}{\mbox{}\hfill $\square$ \\}

\newcommand{\keywords}[1]{{\noindent\bf Key Words }#1}

\ifthenelse{\equal{\mylabelonoff}{on}}
{}{}
\ifthenelse{\equal{\allowdisbrk}{yes}}
{\allowdisplaybreaks}{}

\usepackage[mathscr]{eucal}
\usepackage{color}

\newcommand{\ol}{\overline}

\newcommand{\rz}{\mathbb{R}}

\newcommand{\rt}{\rz^3}
\newcommand{\rttt}{\rz^{3\times3}}
\newcommand{\rttts}{\rttt_{\sym}}

\DeclareMathOperator{\p}{\partial}

\newcommand{\na}{\nabla}

\DeclareMathOperator{\grad}{grad}
\DeclareMathOperator{\Grad}{Grad}

\DeclareMathOperator{\CoGrad}{CoGrad}
\DeclareMathOperator{\curl}{curl}
\DeclareMathOperator{\Curl}{Curl}

\renewcommand{\div}{\operatorname{div}}
\DeclareMathOperator{\Div}{Div}

\newcommand{\trans}[1]{{#1}^t}
\newcommand{\T}{T}
\newcommand{\TS}{S}
\newcommand{\Tp}{P}
\newcommand{\eps}{\varepsilon}

\newcommand{\om}{\Omega}
\newcommand{\dom}{\p\!\om}
\newcommand{\ga}{\Gamma}

\DeclareMathOperator{\sym}{sym}

\DeclareMathOperator{\tr}{tr}
\renewcommand{\skew}{\operatorname{skew}}
\newcommand{\zvec}[2]{\begin{bmatrix}#1\\#2\end{bmatrix}}
\newcommand{\dvec}[3]{\begin{bmatrix}#1\\#2\\#3\end{bmatrix}}
\newcommand{\zmat}[4]{\begin{bmatrix}#1&#2\\#3&#4\end{bmatrix}}

\DeclareMathOperator{\id}{id}
\def\set#1#2{\{#1\,:\,#2\}}

\DeclareMathOperator{\Lebesgue}{\mathsf{L}}
\newcommand{\Lgen}[2]{\Lebesgue^{#1}_{#2}}
\def\Li{\Lgen{\infty}{}}
\def\Liom{\Li(\om)}
\def\Lo{\Lgen{1}{}}
\def\Loom{\Lo(\om)}
\def\Lt{\Lgen{2}{}}
\def\Ltom{\Lt(\om)}

\DeclareMathOperator{\Sobolev}{\mathsf{H}}
\newcommand{\Hgen}[3]{\overset{#3}{\Sobolev}{}^{#1}_{#2}}
\def\Ho{\Hgen{1}{}{}}

\DeclareMathOperator{\Cont}{\mathsf{C}}
\newcommand{\Cgen}[2]{\overset{#2}{\Cont}{}^{#1}}
\def\Ci{\Cgen{\infty}{}}
\def\Cic{\Cgen{\infty}{\circ}}

\def\Coc{\Cgen{1}{\circ}}

\def\Cicom{\Cic(\om)}

\def\Cocom{\Coc(\om)}

\DeclareMathOperator{\dirichlet}{\mathcal{H}}
\newcommand{\qharmdi}[2]{\dirichlet^{#1}_{#2}(\om)}
\newcommand{\harmdi}{\qharmdi{}{}}

\newcommand{\Hggen}[3]{\overset{#2}{\Sobolev}(\grad;#3)}
\newcommand{\HGgen}[3]{\overset{#2}{\Sobolev}(\Grad;#3)}
\newcommand{\Hcgen}[3]{\overset{#2}{\Sobolev}(\curl_{#1};#3)}
\newcommand{\HCgen}[3]{\overset{#2}{\Sobolev}(\Curl_{#1};#3)}
\newcommand{\Hdgen}[3]{\overset{#2}{\Sobolev}(\div_{#1};#3)}
\newcommand{\HDgen}[3]{\overset{#2}{\Sobolev}(\Div_{#1};#3)}
\newcommand{\HCsymgen}[3]{\overset{#2}{\Sobolev}{}_{\sym}(\Curl_{#1};#3)}
\newcommand{\Hgom}{\Hggen{}{}{\om}}
\newcommand{\HGom}{\HGgen{}{}{\om}}
\newcommand{\Hcom}{\Hcgen{}{}{\om}}
\newcommand{\HCom}{\HCgen{}{}{\om}}
\newcommand{\Hdom}{\Hdgen{}{}{\om}}
\newcommand{\HDom}{\HDgen{}{}{\om}}
\newcommand{\Hgcom}{\Hggen{}{\circ}{\om}}
\newcommand{\HGcom}{\HGgen{}{\circ}{\om}}
\newcommand{\Hccom}{\Hcgen{}{\circ}{\om}}
\newcommand{\HCcom}{\HCgen{}{\circ}{\om}}

\newcommand{\Hczom}{\Hcgen{0}{}{\om}}
\newcommand{\HCzom}{\HCgen{0}{}{\om}}
\newcommand{\Hdzom}{\Hdgen{0}{}{\om}}
\newcommand{\HDzom}{\HDgen{0}{}{\om}}
\newcommand{\Hcczom}{\Hcgen{0}{\circ}{\om}}
\newcommand{\HCczom}{\HCgen{0}{\circ}{\om}}

\newcommand{\HCcsymom}{\HCsymgen{}{\circ}{\om}}

\newcommand{\normdst}{\hspace{-0.4ex}}
\newcommand{\scp}[2]{\left\langle#1,#2\right\rangle}

\newcommand{\scpLtom}[2]{\scp{#1}{#2}_{\Ltom}}

\newcommand{\norm}[1]{\left|\normdst\left|#1\right|\normdst\right|}
\newcommand{\dnorm}[1]{\left|\normdst\left|\normdst\left|#1\right|\normdst\right|\normdst\right|}

\newcommand{\normLtom}[1]{\norm{#1}_{\Ltom}}

\newtheorem{lem}{Lemma}

\newtheorem{theo}[lem]{Theorem}
\newtheorem{cor}[lem]{Corollary}
\newtheorem{rem}[lem]{Remark}

\renewcommand{\Tp}{P}  
\renewcommand{\T}{P}   
\renewcommand{\TS}{Q}  
\renewcommand{\trans}[1]{{#1}^T}

\newcommand{\flow}{\Phi}

\begin{document}

\maketitle{}

\begin{abstract}
We prove a Korn-type inequality in $\HCgen{}{\circ}{\om,\rttt}$
for tensor fields $\Tp$ mapping $\om$ to $\rttt$. 
More precisely, let $\om\subset\rt$
be a bounded domain 
with connected Lipschitz boundary $\dom$.
Then, there exists a constant $c>0$ such that 
\begin{equation}
\label{kma}
c\norm{\Tp}_{\Lt(\om,\rttt)}
\leq\norm{\sym\Tp}_{\Lt(\om,\rttt)}+\norm{\Curl\Tp}_{\Lt(\om,\rttt)}
\end{equation}
holds for all tensor fields $\Tp\in\HCgen{}{\circ}{\om,\rttt}$, i.e., 
all $\Tp\in\HCgen{}{}{\om,\rttt}$ with vanishing tangential trace on $\dom$.
Here, rotation and tangential trace are defined row-wise.
For compatible $\Tp$, i.e., $\Tp=\na v$ and thus $\Curl\Tp=0$, 
where $v\in\Ho(\om,\rt)$ are vector fields having components $v_{n}$,
for which $\na v_{n}$ are normal at $\dom$,
the presented estimate \eqref{kma} reduces to a non-standard variant
of Korn's first inequality, i.e.,
$$c\norm{\na v}_{\Lt(\om,\rttt)}\leq\norm{\sym\na v}_{\Lt(\om,\rttt)}.$$
On the other hand, for skew-symmetric $P$, i.e., $\sym\Tp=0$,
\eqref{kma} reduces to a non-standard version of Poincar\'e's estimate. 
Therefore, since \eqref{kma} admits the classical boundary conditions 
our result is a common generalization
of the two classical estimates, namely
Poincar\'e's resp. Korn's first inequality.\\
\end{abstract}


\section{Introduction: Infinitesimal Gradient Plasticity}

The motivation for our new estimate is a formulation 
of infinitesimal gradient plasticity \cite{Reddy06a}. 
Our model is taken from Neff et al. \cite{Neff_Chelminski07_disloc}. 
Let $\om\subset\rt$ be a bounded domain. 
The goal is to find the displacement 
$u:[0,\infty)\times\om\mapsto\rt$ and the possibly non-symmetric plastic 
distortion tensor 
$\Tp:[0,\infty)\times\om\mapsto\rttt$, 
such that in $[0,\infty)\times\om$
\begin{align}
\Div\sigma&=f,
&\sigma&=2\mu\sym(\na u-\Tp)+\lambda\tr(\na u-\Tp)\id,\nonumber\\
\dot{\Tp}&\in\flow(\Sigma),
&\Sigma&=\sigma-2\mu\sym\Tp-\mu L_c^2\Curl\Curl\Tp,\label{infinitesimal_strain_invariant_model}
\intertext{hold. The system is completed by the boundary conditions}
u(t,x)&=0,
&\nu(x)\times\Tp(t,x)&=0\qquad
\forall\,(t,x)\in[0,\infty)\times\dom\nonumber
\end{align}
and the initial condition $\Tp(0,x)=0$ for all $x\in\om$. 
The underlying thermodynamic potential including the plastic gradients
in form of the dislocation density tensor $\Curl\Tp$ is
\begin{align*}
\int_{\om}\mu|&\sym(\na u-\Tp)|^2
+\frac{\lambda}{2}|\tr(\na u-\Tp)|^2
-f\cdot u+\mu|\sym\Tp|^2
+\frac{\mu}{2}L_c^2|\Curl\Tp|^2.
\end{align*}

Here, $\mu,\lambda$ are the elastic Lam\'e moduli and $\sigma$ 
is the symmetric Cauchy stress tensor. 
The system is driven by nonzero body forces denoted by $f$.
The exterior normal to the boundary $\dom$ is denoted by $\nu$ 
and the plastic distortion $\Tp$ is required to satisfy row-wise 
the homogeneous tangential boundary condition 
which means that the boundary $\dom$ is a 
perfect conductor regarding the plastic distortion.
\footnote{This homogeneous tangential boundary condition on $\Tp$ 
is consistent with $\nu\times\na u=0$ on $\dom$ 
which follows from $u=0$ on $\dom$.}

Moreover, $\flow:\rttt\mapsto\rttt$ is the monotone,
multivalued flow-function with $\flow(0)=0$ 
and $\flow(\rttts)\subset\rttts$.
In general, $\Sigma$ is not symmetric even if $\Tp$ is symmetric. 
Thus, the plastic inhomogeneity is responsible for 
the plastic spin (the possible non-symmetry of $\Tp$). 
The mathematically suitable space 
for symmetric plastic distortion $\Tp$ is the classical space 
$\Hcom$ for each row of $\Tp$ \cite{Neff_Sydow_Wieners08,Reddy06a}. 
This case appears when choosing $\flow:\rttt\mapsto\rttts$.

In the large scale limit $L_c\to0$ we recover 
a classical elasto-plasticity model with local kinematic hardening 
and symmetric plastic strain $\eps_p:=\sym\Tp$, 
since then $\dot{\Tp}\in\rttts$. 

Uniqueness of classical solutions for rate-independent 
and rate-dependent formulations of this model is shown in \cite{Neff_Chelminski07_disloc}. 
The more difficult existence question for the rate-independent model 
in terms of a weak reformulation is addressed in \cite{Neff_Chelminski07_disloc}. 
First numerical results for a simplified rate-independent irrotational formulation 
(no plastic spin, i.e., symmetric plastic distortion $\Tp$) 
are presented in \cite{Neff_Sydow_Wieners08}, cf \cite{Reddy06}.  
In \cite{Ebobisse_Neff09} the model has been extended 
to rate-independent isotropic hardening based on the concept 
of a dissipation function defined in terms of the equivalent plastic strain. 
From a modeling point of view, 
it is strongly preferable to again have only the symmetric (rate) part 
of the plastic distortion $P$ appear in the dissipation potential.

The existence and uniqueness can be settled 
by recasting the model as a variational inequality, 
if it is possible to define a bilinear form 
which is coercive with respect to appropriate spaces. 
This program has been achieved for other variants 
of the model in \cite{Ebobisse_Neff09}. 
It had to remain basically open 
for the above system \eqref{infinitesimal_strain_invariant_model}. 
In this case, the appropriate space for the plastic distortion $\Tp$ is the completion 
$\HCcsymom$ of the linear space 
$$\set{\Tp\in\Ci(\ol{\om},\rttt)}{\Tp_{n}\text{ normal at }\dom,\,n=1,2,3}$$
with respect to the norm $\dnorm{\,\cdot\,}$, 
where $\Tp_{n}$ are the columns of $\trans{\Tp}$ and
$$\dnorm{\Tp}^2:=\normLtom{\sym\Tp}^2+\normLtom{\Curl\Tp}^2.$$
Despite first appearance, 
this quadratic form indeed defines a norm 
as shown in \cite{Neff_Chelminski07_disloc}.
Thus $\HCcsymom$ is a Hilbert-space. 
However, in this space it is not immediately obvious 
how to define a linear and bounded tangential trace operator. 
Since only $\normLtom{\sym\Tp}$ appears, 
it is also not clear, how to control the skew-symmetric part of $\Tp$. 
Therefore, the crucial embedding 
$$\HCcsymom\subset\Ltom$$ 
is not clear as well. 
As a consequence of our main result of this paper we obtain that nevertheless
$$\HCcsymom=\HCcom$$
holds with equivalent norms in case the domain $\om$ 
is simply connected and has a Lipschitz boundary. 
The result of this paper 
has been announced in \cite{neffpaulywitschgenkornpamm}.

For the proof of our main result \eqref{kma} 
we combine techniques from electro-magnetic and elastic theory,
namely the Helmholtz decomposition, the Maxwell compactness property
and Korn's inequality. 
Their basic variants are well known results
which can be found in many books, e.g., 
in \cite{leisbook} and the literature cited there.
More sophisticated and related versions are presented, e.g.,
in \cite{paulydeco,picardcomimb,picarddeco,picardweckwitschxmas,weckmax} 
for Maxwell's equations and \cite{Ciarlet10,Neff00b} for Korn's inequality. 

This paper is organized as follows. 
After this motivation we introduce our notation, 
definitions and provide some background results.
In section \ref{mainresultsec} we give the proof for our main estimates. 
In the last section \ref{twodimsec} 
we establish a connection to a related result by Garroni et al. \cite{Garroni10} 
for the two-dimensional case.

\section{Definitions and Preliminaries}

Let $\om$ be a bounded domain in $\rt$
with connected Lipschitz continuous boundary $\ga:=\dom$.

\subsection{Functions and Vector Fields}

The usual Lebesgue spaces of square integrable functions, 
vector or tensor fields on $\om$ 
with values in $\rz$, $\rt$ or $\rttt$, respectively,
will be denoted by $\Ltom$.
Moreover, we introduce the standard Sobolev spaces 
\begin{align*}
\Hgom&=\set{u\in\Ltom}{\grad u\in\Ltom},
&\norm{u}_{\HGom}^2&:=\normLtom{u}^2+\normLtom{\grad u}^2,\\
\Hcom&=\set{v\in\Ltom}{\curl v\in\Ltom},
&\norm{v}_{\Hcom}^2&:=\normLtom{v}^2+\normLtom{\curl v}^2,\\
\Hdom&=\set{v\in\Ltom}{\div v\in\Ltom},
&\norm{v}_{\Hdom}^2&:=\normLtom{v}^2+\normLtom{\div v}^2.
\end{align*}
$\Hgom$ is often denoted by $\Hgen{1}{}{}(\om)$.
Furthermore, we define their closed subspaces
$\Hgcom$, $\Hccom$ as completition under the respective norms of
the scalar resp. vector valued space $\Cicom$
of compactly supported and smooth test functions resp. vector fields.
In the latter Sobolev spaces
the usual homogeneous scalar resp. tangential boundary conditions
$$u|_{\ga}=0,\quad\nu\times v|_{\ga}=0$$
are generalized, where $\nu$ denotes the outer unit normal at $\ga$. 
We note in passing that $\nu\times v|_{\ga}=0$ is equivalent to
$\tau\cdot v|_{\ga}=0$ for all tangential directions $\tau$ at $\ga$,
which means that $v$ is normal at $\ga$. 
Furthermore, we need the spaces of irrotational or solenoidal vector fields
\begin{align*}
\Hczom&:=\set{v\in\Hcom}{\curl v=0},\\
\Hcczom&:=\set{v\in\Hccom}{\curl v=0}, \\
\Hdzom&:=\set{v\in\Hdom}{\div v=0},
\end{align*}
where the index $0$ indicates vanishing $\curl$ or $\div$, respectively.
All these spaces are Hilbert spaces. 
E.g., in classical terms we have $v\in\Hcczom$, if and only if
$$\curl v=0,\quad\nu\times v|_{\ga}=0.$$
For an introduction of these spaces see 
\cite[p. 11-12, 148]{leisbook} or \cite[p. 26]{Raviart79}. 
The most important tool for our analysis is the compact embedding
\begin{align*}
\Hccom\cap\Hdom&\hookrightarrow\Ltom,
\end{align*}
which is often referred as `Maxwell compactness property', 
see \cite[p. 158]{leisbook} and
\cite{weckmax,picardcomimb, webercompmax,witschremmax,picardweckwitschxmas}. 
A first immediate consequence is that the space 
of so called `harmonic Dirichlet fields'
\begin{align*}
\harmdi:=\Hcczom\cap\Hdzom
\end{align*}
is finite dimensional. 
A vector field $v$ belonging to $\harmdi$ means in classical terms that 
$$\curl v=0,\quad\div v=0,\quad\nu\times v|_{\ga}=0.$$
The dimension of $\harmdi$ equals the second Betti number of $\om$,
see \cite[p. 159]{leisbook} and \cite[Theorem 1]{picardboundaryelectro}.
Since we assume the boundary $\ga$ to be connected,
there are no Dirichlet fields besides zero, i.e., 
$$\harmdi=\{0\}.$$
This condition on the domain $\om$ resp. its boundary $\ga$
is satisfied e.g. for a ball or a torus.

By a usual indirect argument
we achieve another immediate consequence, 
see \cite[p. 158, Theorem 8.9]{leisbook} or \cite[Lemma 3.4]{Raviart79}: 

\begin{lem}
\label{poincaremax}
{\sf(Maxwell Estimate for Vector Fields)}
There exists a positive constant $c_{m}$, 
such that for all vector fields $v\in\Hccom\cap\Hdom$
$$\normLtom{v}\leq 
c_{m}\big(\normLtom{\curl v}^2+\normLtom{\div v}^2\big)^{1/2}.$$
\end{lem}

By definition of the weak divergence, 
the projection theorem and Rellich's selection theorem 
\cite[p. 14]{leisbook} we have from 
\cite[p. 148, Theorem 8.3]{leisbook} 
or \cite[Lemma 3.5]{weberregmax}, \cite[Theorem 3.45]{monkbook}

\begin{lem}
\label{helmdeco}
{\sf(Helmholtz Decomposition for Vector Fields)}
We have the orthogonal decomposition
$$\Ltom=\grad\Hgcom\oplus\Hdzom.$$
\end{lem}
 
\subsection{Tensor Fields}

We extend our calculus to $(3\times3)$-tensor (matrix) fields.
For vector fields $v$ with components in $\Hgom$
and tensor fields $\T$ with rows in $\Hcom$ resp. $\Hdom$, i.e.,
$$v=\dvec{v_{1}}{v_{2}}{v_{3}},\quad v_{n}\in\Hgom,\quad
\trans{\T}=[\T_{1}\,\T_{2}\,\T_{3}],\quad
\T_{n}\in\Hcom\text{ resp. }\Hdom$$
we define
$$\Grad v:=\dvec{\trans{\grad}v_{1}}{\trans{\grad}v_{2}}{\trans{\grad}v_{3}}=J_{v}=\na v,\quad
\Curl\T:=\dvec{\trans{\curl}\T_{1}}{\trans{\curl}\T_{2}}{\trans{\curl}\T_{3}},\quad
\Div\T:=\dvec{\div\T_{1}}{\div\T_{2}}{\div\T_{3}},$$ 
where $J_{v}$ denotes the Jacobian of $v$ 
and $\trans{}$ the transpose.
We note that $v$ and $\Div\T$ are vector fields, whereas
$\T$, $\Curl\T$ and $\Grad v$ are tensor fields.
The corresponding Sobolev spaces will be denoted by 
$\HGom$, $\HGcom$, $\HCom$, $\HCcom$, $\HCzom$, $\HCczom$, $\HDom$, $\HDzom$.
As usual, we denote by $\sym\T:=1/2(\T+\trans{\T})$
the symmetric part of a tensor $\T$.

Let us now present our three crucial tools to prove the new estimate.
First we have obvious consequences 
from Lemmas \ref{poincaremax} and \ref{helmdeco}:

\begin{cor}
\label{poincaremaxten}
{\sf(Maxwell Estimate for Tensor Fields)}
For all $\T\in\HCcom\cap\HDom$ 
$$\normLtom{\T}\leq 
c_{m}\big(\normLtom{\Curl\T}^2+\normLtom{\Div\T}^2\big)^{1/2}.$$
\end{cor}

\begin{cor}
\label{helmdecoten}
{\sf(Helmholtz Decomposition for Tensor Fields)}
We have the orthogonal decomposition
$$\Ltom=\Grad\HGcom\oplus\HDzom.$$
\end{cor}

The third important tool is Korn's first inequality
\cite[p. 207]{leisbook} or \cite[p. 54]{Valent88}:

\begin{lem}
\label{korn}
{\sf(Korn's First Inequality)}
For all $v\in\HGcom$
$$\normLtom{\Grad v}\leq\sqrt{2}\normLtom{\sym\Grad v}.$$
\end{lem}

\section{Main Results}
\label{mainresultsec}

For tensor fields $\T\in\HCom$ we define the semi-norm
\begin{align*}
\dnorm{\T}:=\big(\normLtom{\sym\T}^2+\normLtom{\Curl\T}^2\big)^{1/2}.
\end{align*}

\begin{lem}
\label{mainlem}
Let $\hat{c}:=\max\{2,\sqrt{5}c_{m}\}$. Then, for all $\T\in\HCcom$ 
$$\normLtom{\T}\leq\hat{c}\dnorm{\T}.$$
\end{lem}

\noindent{\bf Proof }
Let $\T\in\HCcom$.
According to Corollary \ref{helmdecoten} we orthogonally decompose 
$$\T=\Grad v+\TS\in\Grad\HGcom\oplus\HDzom.$$
Then, $\Curl\T=\Curl\TS$ and we observe $\TS\in\HCcom\cap\HDzom$ since 
\begin{align}
\label{gradsubsetcurl}
\Grad\HGcom\subset\HCczom.
\end{align}
By Corollary \ref{poincaremaxten} we have
\begin{align}
\label{estpsi}
\normLtom{\TS}\leq c_{m}\normLtom{\Curl\T}.
\end{align}
Then, by Lemma \ref{korn} and \eqref{estpsi} we obtain easily
\begin{align*}
\normLtom{\T}^2
&=\normLtom{\Grad v+\TS}^2
=\normLtom{\Grad v}^2+\normLtom{\TS}^2\\
&\leq2\normLtom{\sym\Grad v}^2+\normLtom{\TS}^2
=2\normLtom{\sym(\T-\TS)}^2+\normLtom{\TS}^2\\
&\leq4\normLtom{\sym\T}^2+5\normLtom{\TS}^2
\leq4\normLtom{\sym\T}^2+5c_m^2\normLtom{\Curl\T}^2.
\qquad\Box
\end{align*}

The immediate consequence is

\begin{theo}
\label{maintheo}
On $\HCcom$ the norms $\norm{\,\cdot\,}_{\HCom}$
and $\dnorm{\,\cdot\,}$ are equivalent.
In particular, $\dnorm{\,\cdot\,}$ is a norm on $\HCcom$ and 
$$\exists\,c>0\quad\forall\,\T\in\HCcom\qquad
c\norm{\T}_{\HCom}\leq\normLtom{\sym\T}+\normLtom{\Curl\T}.$$
\end{theo}

Setting $\T:=\Grad v$ we obtain 
by Lemma \ref{mainlem} and \eqref{gradsubsetcurl}

\begin{rem}
\label{genkorn}
{\sf(Korn's First Inequality: Tangential-Variant)}
For all $v\in\HGcom$ 
\begin{align}
\label{kornineq}
\normLtom{\Grad v}\leq\hat{c}\normLtom{\sym\Grad v}.
\end{align}
This is Korn's first inequality from Lemma \ref{korn}
with a larger constant $\hat{c}$.
Since the boundary $\ga$ is connected, i.e.,
$\harmdi=\{0\}$, we have
$\Grad\HGcom=\HCczom$. Thus,
\eqref{kornineq} holds for all
$v\in\HGom$ with $\Grad v\in\HCczom$, i.e.,
with $\Grad v_{n}$, $n=1,2,3$, normal at $\ga$,
which then extends Lemma \ref{korn} 
through the (apparently) weaker boundary condition.
\end{rem}

\section{Two-Dimensions: a Result of Garroni et al.}
\label{twodimsec}

Let $\om$ be a bounded domain in $\rz^2$
with connected Lipschitz continuous boundary $\ga$,
which is equivalent (in $\rz^2$) to the topological property
that $\om$ is simply connected.
For tensor fields $\T:\om\mapsto\rz^{2\times2}$ 
we define analogously the $\Curl$-operator by
$$\Curl\T
=\Curl\zmat{\T_{11}}{\T_{12}}{\T_{21}}{\T_{22}}
=\zvec{\curl\trans{[\T_{11}\,\T_{12}]}}{\curl\trans{[\T_{21}\,\T_{22}]}}
=\zvec{\p_{1}\T_{12}-\p_{2}\T_{11}}{\p_{1}\T_{22}-\p_{2}\T_{21}},$$
where now $\curl$ denotes the two dimensional scalar rotation
and $\Curl\T$ is a vector.
With the appropriate changes, 
Lemma \ref{mainlem} and Theorem \ref{maintheo} hold as well.
In particular, there exists a positive constant $c$, such that
$$c\normLtom{\T}\leq\normLtom{\sym\T}+\normLtom{\Curl\T}$$
holds for all $\T\in\HCcom$. 

During the preparation of our paper 
we got aware that a two-dimensional related result 
may be inferred from Garroni et al. \cite{Garroni10}. 
Instead of tangential boundary conditions $\nu\times P|_{\ga}=0$ 
they impose the normalization condition 
\begin{align}
\label{intskew}
\int_\om\skew\T=0.
\end{align}
Let us define the total variation measure 
of the distribution $\Curl\T$ for $\T\in\Loom$ by
\begin{align*}
|\Curl\T|_{\om}
:=\sup_{\substack{v\in\Cocom\\\norm{v}_{\Liom}\leq1}}\scpLtom{\T}{\CoGrad v},\quad   
\CoGrad v
:=\zmat{\p_{2}v_{1}}{-\p_{1}v_{1}}{\p_{2}v_{2}}{-\p_{1}v_{2}}.
\end{align*}
We note
$$\scpLtom{\T}{\CoGrad v}
=\int_{\om}\T_{11}\p_{2}v_{1}-\T_{12}\p_{1}v_{1}
+\T_{21}\p_{2}v_{2}-\T_{22}\p_{1}v_{2}.$$
Using partial integration, i.e., 
$\scpLtom{\T}{\CoGrad v}=\scpLtom{\Curl\T}{v}$ for $v\in\Cocom$,
it is easy to see that $|\Curl\T|_{\om}=\norm{\Curl\T}_{\Loom}$ if $\Curl\T\in\Loom$. 
In \cite[Theorem 9]{Garroni10} it is shown 
that for $\om$ having a Lipschitz boundary 
and a special `slicing' property, 
there exists a constant $c>0$, such that
$$c\normLtom{\T}\leq\normLtom{\sym\T}+|\Curl\T|_{\om}$$
holds for all $\T\in\Loom$ with \eqref{intskew}.
Their proof uses essentially 
that in $\rz^2$ the operators $\curl$ and $\div$
can be exchanged by the simple transformation, i.e., 
$\curl\trans{[v_1,v_2]}=\div\trans{[-v_2, v_1]}$.
Thus, such a strong result may not be true in higher space dimensions $N\geq3$ 
and it is open whether the normalization condition 
\eqref{intskew} can be exchanged with the more natural 
tangential boundary conditions.

{\footnotesize
\bibliographystyle{plain} 
\bibliography{/Users/paule/Library/texmf/tex/TeXinput/bibtex/paule,/Users/paule/Library/texmf/tex/TeXinput/bibtex/literatur1}
\vspace*{5mm}
\begin{center}
{\sf\begin{tabular}{l}
Patrizio Neff, Dirk Pauly, Karl-Josef Witsch\\ 
Universit\"at Duisburg-Essen, Fakult\"at f\"ur Mathematik, Campus Essen\\
Universit\"atsstr. 2, 45117 Essen, Germany\\ 
e-mail: patrizio.neff@uni-due.de, dirk.pauly@uni-due.de, kj.witsch@uni-due.de
\end{tabular}}
\end{center}
}

\end{document}